\documentclass[12pt,leqno]{article}

\begin{document}

\title{Complexity and ordinary life, \\ and some mathematics in both}

\author{Stephen Semmes   \\
	Department of Mathematics \\
	Rice University}

\date{}

\maketitle

	What is mathematics, exactly?  This is a somewhat complicated
question, with no simple answer.  In any event, mathematics is like a
large place, with many regions and villages, and many different ways
of doing things.  One can also try to make up new ways of doing
things, in connection with whatever might be of interest.

	Let us look here at a few points which can come up naturally
in ordinary life, as well as involve some substantial mathematics (in
some of their forms).

\subsection* {Differences between answers of ``yes'' and ``no''}

	Imagine the following kind of question: ``Does so-and-so have
a pencil in his or her office?''

	If someone finds a pencil in the office, then that provides a
way in which an answer of ``yes'' can be clearly established.  The pencil 
can simply be shown.  An answer of ``no'' is quite different, and apparently
more complicated.  How can one establish an answer of ``no'', without just
going through all of the contents of the office?

	This type of distinction is sometimes described technically
through phrases like \emph{effective witnesses} or \emph{succinct
certificates}.  Roughly speaking, an effective witness or succinct
certificate is something which would be used to establish an answer to
a given question, which is reasonably small or manageable, and for
which the verification of this as providing an answer would be fairly
easy and definite.

	This general idea can be given precise forms in suitable
mathematical contexts, and in particular it has a basic role in
theoretical computer science.  (Compare with
\cite{aho-hopcroft-ullman, hopcroft-ullman, lewis-papadimitriou,
papadimitriou}.)

	In the example above, an actual pencil in the office would
serve as an effective witness or succinct certificate for an answer of
``yes'', to the question of whether there is a pencil in so-and-so's
office.  For an answer of ``no'', it is not clear what might make
sense as an effective witness.  There may not be one.

	This is a very general issue, and one that comes up in many
ways.  Instead of a pencil, the question might be more interesting,
concerning the possibility that so-and-so is in possession of items of
value to sports fans, or nice jewelry, or whether so-and-so remembers
about Wednesday evening.  One might really want to be able to have a
definite answer to the question, and that might not be so easy to come
by.  (So-and-so might thus be in a good position to be evasive.)

	Even if there are effective witnesses for a given answer
(``yes'' or ``no''), this is not the same as saying that it is easy to
find an effective witness.  In other words, a pencil or other object
might serve as an effective witness for an answer of ``yes'' if one
can find one, but that does not mean that it is easy to do so, or that
one has a good method for doing so.  One might be faced once more with
the prospect of something like an exhaustive search, and that might
not be appealing or feasible.  (Again, this could be helpful for
so-and-so.)

	In an argument, this might come up in a slightly different
(but equivalent) way, in which the roles of ``yes'' and ``no'' are
reversed.  One person might want to put forward a general statement,
while another person might want to look in to its correctness.  The
statement could be something like ``Members of the Oblidian Club are
all sipifsts!''  An effective witness for an answer of ``no'' could be
a person in the Oblidian Club who is manifestly \emph{not} a sipifst.
It may not be easy to locate such an individual even if he or she
exists, but at least the possibility is there.  For an answer of
``yes'', there might not be any simple effective witness like this at
all.

	These are phenomena which occur all the time, in arguments
between people in particular.  In theoretical computer science, there
are mathematically-precise versions of these notions, and in fact
there are famous unsolved problems pertaining to them.  In technical
terms, one of these is the problem that asks whether the ``complexity
class'' NP is equal to the complementary class co-NP.  The NP class
involves questions for which effective witnesses for answers of
``yes'' always exists, while the co-NP class entails questions for
which effective witnesses for ``no'' exist.  

	Even if effective witnesses exist, it might not be easy to
find them, and the complexity class P consists of questions for which
there is a method to figure out an answer of ``yes'' or ``no'' in a
limited amount of time (namely, in ``polynomial time'').  It is a
famous unsolved problem to know whether the class P might actually be
the same as NP, even though P seems to be significantly more
restrictive than NP.  

	See \cite{aho-hopcroft-ullman, hopcroft-ullman,
lewis-papadimitriou, papadimitriou} for more about these classes and
unsolved problems.

	At any rate, versions of these issues come up a lot, and in
various forms, in both ordinary activities, and in more technical or
mathematical situations.  We shall see a bit more of this later.

\subsection* {Hamburgers and exponentiation}

	Some years ago, there was a restaurant (with special emphasis on
hamburgers) which had a sign that said something like the following:
\begin{quote}
	We have $256$ different kinds of hamburgers.
\end{quote}
This might seem to be a rather remarkable statement, but in fact we
can look at it mathematically and see how it makes a lot of sense.

	The first point is to understand something about what the
number $256$ actually is.  It is not just \emph{any} number, but a
very special one.

	One might first notice that $256$ is an even number, meaning
that it is divisible by $2$.  One can see this quickly because its
last digit is an even number.  (One is told in school, at least in
years gone by, that a number is even if its last digit is even.  This
is not hard to verify anyway.)  Thus $256$ can be written as $2$ times
something.  One can check that that something is $128$:
\begin{displaymath}
	256 = 2 \times 128.
\end{displaymath}

	This is a bit curious, because $128$ is also an even number.  In 
fact $128$ is $2$ times $64$, and so we get that
\begin{displaymath}
	256 = 2 \times 2 \times 64.
\end{displaymath}

	Now again, $64$ is even.  It is $2$ times $32$.  And $32$ is even,
$32 = 2 \times 16$.  And then $16 = 2 \times 8$, $8 = 2 \times 4$, and 
$4 = 2 \times 2$.  

	In the end, we get that $256$ can be obtained by multiplying
together a bunch of $2$'s, without needing any other numbers.
Specifically, $256$ is the product of eight $2$'s, i.e.,
\begin{displaymath}
    256 = 2 \times 2 \times 2 \times 2 \times 2 \times 2 \times 2 \times 2.
\end{displaymath}
This is sometimes written as 
\begin{displaymath}
	256 = 2^8.
\end{displaymath}
In general, $2^m$ means the number that one gets by multiplying together
$m$ $2$'s.  This process is called \emph{exponentiation}, and the number
$m$ in $2^m$ is called the \emph{exponent}.

	This is analogous to the relationship between multiplication
and addition.  That is, $m \times 2$ is the same as \emph{adding}
together $m$ $2$'s.  With $2^m$, one \emph{multiplies} the $2$'s
instead of adding them.

	One of the interesting features of exponentiation is that it
leads to pretty big numbers rather quickly.  We shall see more of this
soon.  For the moment, let us continue with the question of the
hamburgers.  Now that we know that $256$ is the same as eight $2$'s
multiplied together, what does that suggest about the restaurant
having $256$ different kinds of hamburgers?

	In fact, there is a rather simple answer to this.  Consider the
following $8$ choices that one might be offered, in having a hamburger
at the restaurant:
\begin{enumerate}

\item	Do you want pickles?

\item 	Do you want lettuce?

\item 	Do you want tomatoes?

\item	Do you want onions?

\item	Do you want cheese?

\item	Do you want ketchup?

\item	Do you want mustard?

\item	Do you want mayonnaise?

\end{enumerate}
Each of these choices represents $2$ possibilities, ``yes'' or ``no''.
With all $8$ choices, where one is free to give answers of ``yes'' or
``no'' for each one, independently of the rest, one obtains a total of
$256$ possible types of hamburgers to have.  This comes from the fact
that $256$ is really the same as multiplying eight $2$'s together,
where each of the eight $2$'s corresponds to the $2$ options for one
of the questions above.  (This is a standard observation, and it is
not too hard to check.)

	This is not to say that the restaurant maintains a supply of
all $256$ different types of hamburgers at any given moment.  Instead,
the restaurant has a supply of each of the $8$ different ingredients
mentioned in the questions above.  They might add the ingredients in
preparing the hamburger to be served, or they might have a place where
a customer can get his or her own ingredients.

	One might say that the restaurant has $256$ different kinds of
hamburgers \emph{implicitly}.  This is as opposed to having them all
\emph{explicitly}, with $256$ actual hamburgers present, covering all
$256$ possible different ways of answering the $8$ questions above.

\subsection* {Hamburgers and exponentiation, part 2}

	Now let us imagine that there might be more than the $8$
choices mentioned above.  There might be options for oregano, bean
sprouts, or mushrooms, for instance, or a whole-wheat bun.  There
might also be an option for a vegetarian burger instead of a usual
hamburger.

	Suppose that there are $20$ yes-or-no choices for the burgers,
rather than $8$.  This means that there would be $2^{20}$ different
kinds of burgers, i.e., a product of twenty $2$'s.  This number is
equal to
\begin{displaymath}
	1,048,576
\end{displaymath}
(i.e., a little more than a million).

	What if there were $30$ such choices instead?  This would mean
that there would be $2^{30}$ different kinds of burgers.  This number is
equal to
\begin{displaymath}
	1,073,741,824
\end{displaymath}
(a little more than a billion).

	Each time that one adds ten more choices, one should multiply
the total number of burgers by $2^{10}$, which is $1024$ (a little
more than a thousand).  With $40$ choices one would get about a
trillion different kinds of hamburgers, and with $50$ one would get
about a quadrillion.

	Relatively speaking, it would not be that hard to make a place
where one really could have $20$, $30$, or more choices for what to
have on the burger.  It might be a little difficult to put in all $30$
types of ingredients on a single hamburger, but still, one can imagine
things like this.  (There are modest variants of this scenario which
might be easier for making the hamburgers, and which still lead to
large numbers of different types.)

	With larger numbers of choices like this, the
\emph{implicitness} becomes more and more of an issue.  One can
imagine a place where a large number of options for the different
hamburgers is possible, so that very many different kinds of
hamburgers are possible, but this is very different from actually
having all of those millions or billions of hamburgers right there.
The range of all of these hamburgers would be represented
\emph{implicitly} in this way, if not \emph{explicitly}.

	To put this into perspective, how many seconds are there in a
year?  There are $60$ seconds in a minute, $60$ minutes in an hour,
$24$ hours in a day, and $365$ days in a year.  Thus the total number
of seconds in a year is
\begin{displaymath}
	60 \times 60 \times 24 \times 365 = 31,536,000.
\end{displaymath}
One might also think about issues of storage for the millions or billions
of hamburgers.  Having a facility in which one can make hamburgers with
any of the various choices is much easier by comparison.

	All of this makes sense much more generally.  A modest number
of independent choices or options can lead to huge numbers of total
different outcomes, like different possible hamburgers, as above.
This works just as well for other types of objects, in addition to
hamburgers.

	As another basic scenario, imagine having a bunch of places in
some area, like towns and villages, and some roads connecting them.
One can then have a lot of different paths that one could take, following
roads between the different places.  I.e., does one go here first and then
there, or the other way around, and so on.

	In general situations like this, the number of different paths
that would be possible can grow exponentially, as in the case of
hamburgers.  Now it would grow exponentially as compared to the number
of places and roads in the region, rather than the number of
ingredients, like pickles, cheese, and mustard.  With a modest number
of places and roads, there can be an enormous number of different
paths that one could take, just as before.

	For this statement, we might agree to only count paths that do
not go in circles anywhere.  Otherwise, one might simply go around
circles over and over again, which leads to numerous paths by itself.
Without ever going around in circles in any given path, one can still
have very many paths, with exponential size for the number of paths.
If one includes paths that can go around in circles (for a while,
say), then one gets even more.

	In this scenario the idea of \emph{implicitness} comes up too.
One might depict the places and roads on a map or a diagram, and this
implicitly gives a way to represent all of the paths on them.  Any
given path can be traced out on the map or diagram, for instance.  An
actual listing of all of the individual paths would be quite
different.  In general there could be too many of them to list in a
reasonable or practical way.

\subsection* {Big sets and effective witnesses}

	Once one has very big sets like this, as in the illustrations
of hamburgers and paths, one has again issues like those from before,
with effective witnesses, complexity, and so on.

	A classical example of this occurs with the case of paths.
Imagine that one has a list of places in the region, and that one asks
the following type of question:
\begin{quote}
	Is there a path along these roads, which visits all of the locations
       	on this list, but where the total distance travelled is at most
       	$230$ miles?
\end{quote} 
For an answer of ``yes'', there is a clear effective witness for this
problem, namely a path with the given properties.  In particular, for
a given path, it is not too hard to check whether or not it has all of
the required features.

	However, this is not to say that it is easy to find such a
path when it exists, or to determine whether such a path does exist
(without necessarily exhibiting it).  Also, it is not clear that there
should be reasonable effective witnesses for an answer of ``no''.  In
other words, how might one be convinced that such a path does not
exist, without simply checking all of the individual paths of total
length at most $230$ miles?

	The problem of finding a yes-or-no answer to the question
above is sometimes called ``the travelling salesperson problem''.  It
has been much studied, and it is a famous example of a problem which
is ``complete'' for the class NP that was mentioned earlier.  See
\cite{aho-hopcroft-ullman, hopcroft-ullman, lewis-papadimitriou,
papadimitriou}.

	There are similar issues that one can see in the context of
hamburgers.  One could have a list of conditions that one would like a
hamburger to satisfy (in terms of the different ingredients that are
available, as before), and one might ask whether such a hamburger
exists.  For any given hamburger, it could be easy to decide whether
that one enjoys the right properties, but this is very different from
being able to say much about what happens among all hamburgers.  In
particular, it may not be clear how to find a hamburger with the given
features, except for just searching through them directly.  There may
not be any good way to know when there is no such hamburger, without
searching through them all.  (When there is a suitable hamburger, one
can stop searching as soon as one finds it, but when there is not, one
may not know when to stop, without going to the end.)

	This type of question about hamburgers is a version of the
``satisfiability'' problem, discussed in \cite{aho-hopcroft-ullman,
hopcroft-ullman, lewis-papadimitriou, papadimitriou}.

	There are a lot of questions and phenomena along these lines.
In technical versions, one can see this in \cite{aho-hopcroft-ullman,
carbone-semmes, hopcroft-ullman, lewis-papadimitriou, papadimitriou},
for instance.  There are also plenty of forms of this which come up in
ordinary life.

	In some ways this is not really a big deal.  I.e., it is just
fair, and reasonable.  On the other hand, it is basic and substantial,
and can be quite tricky.

\end{document}